 \documentclass[final,1p,times]{elsarticle}

\usepackage{amssymb}
\newtheorem{thm}{Theorem}
\newtheorem{theorem}[thm]{Theorem}

\newtheorem{definition}[thm]{Definition}

\newtheorem{lemma}[thm]{Lemma}
\newtheorem{corollary}[thm]{Corollary}
\newcommand{\be}{\begin{equation}}
\newcommand{\ee}{\end{equation}}
\newcommand{\bea}{\begin{eqnarray}}
\newcommand{\eea}{\end{eqnarray}}
\newcommand{\bd}{\begin{displaymath}}
\newcommand{\ed}{\end{displaymath}}

\newcommand{\op}{ \oplus_q  }
\newcommand{\om}{ \ominus_q  }

\newcommand{\lb }{ \left( }
\newcommand{\rb }{ \right ) }

\journal{Journal of Computational and Applied Mathematics }

\begin{document}

\begin{frontmatter}


\author{Mahouton Norbert Hounkonnou\corref{cor1}${}^1$}
\cortext[cor1]{norbert.hounkonnou@cipma.uac.bj}
 \author{Sama Arjika\corref{cor2}${}^1$}
\cortext[cor2]{rjksama2008@gmail.com}
\author{ Won Sang Chung\corref{cor3}${}^2$ }
\cortext[cor3]{mimip4444@hanmail.net}


\title{\bf New families of 
$q$ and $(q;p)-$Hermite  polynomials 
  }
\address{${}^1$International Chair of Mathematical Physics
and Applications \\
 (ICMPA-UNESCO Chair), University of
Abomey-Calavi,\\
 072 B. P.: 50 Cotonou, Republic of Benin,\\
${}^2$Department of Physics and Research Institute of Natural Science, \\
College of Natural Science, \\
Gyeongsang National University, Jinju 660-701, Korea
}

\begin{abstract}
In this paper, we construct 
a new family of $q-$Hermite polynomials  denoted by
$H_n(x,s|q).$ Main properties and  relations 
are 
established and proved. 
In addition, is deduced a sequence of novel polynomials, 
$\mathcal{L}_n(\cdot,\cdot |q),$ which appear to be  connected with well known $(q,n)-$exponential functions $E_{q,n}(\cdot)$\, introduced by Ernst in his work entitled: {\it A New Method for $q-$calculus,} (Uppsala Dissertations in 
Mathematics, Vol. {\bf 25}, 2002). Relevant  results spread in the literature 
are retrieved  as particular cases. Fourier integral transforms  are explicitly computed and discussed. A $(q;p)-$extension of the $H_n(x,s|q)$ is also provided.
\end{abstract}

\end{frontmatter}
\section{Introduction}
\label{int1}
The classical orthogonal polynomials   and the quantum orthogonal polynomials, 
also called $q-$orthogonal polynomials, constitute an interesting set of 
special functions.
  Each family  of   these polynomials occupies different levels
within the so-called Askey-Wilson scheme (Askey and Wilson, 1985; Koekoek and Swarttouw, 1998; 
Lesky, 2005; Koekoek et {\it al}, 2010).
In this scheme, the   Hermite polynomials
$H_n(x)$ are the ground level and are characterized 
by a set of properties: 
(i) they are solutions of   a hypergeometric second order differential equation, 
(ii) they are generated by  a recursion relation,
(iii) they are orthogonal with  respect to  a weight function and
(iv) they obey  the Rodrigues-type formula.
Therefore, there are many ways to construct the Hermite
polynomials.  However, they are more commonly  
deduced from their
  generating function, i.e.,
\be
\sum_{n=0}^\infty\frac{{\bf H}_n(x)}{n!}t^n=e^{2xt-t^2}
\ee
giving rise to the so-called   {\it physicists} Hermite polynomials  \cite{Johann Cigler}.
Another family of Hermite polynomials, called 
the {\it probabilists} Hermite polynomials, is defined as  \cite{Johann Cigler}
 \be
\sum_{n=0}^\infty\frac{H_n(x)}{n!}t^n=e^{xt-\frac{t^2}{2}}.
\ee

The Hermite polynomials are at the bottom of a 
large class of  hypergeometric polynomials
to which most of their properties can be 
generalized \cite{Bukweli&Hounkonnou12a}, \cite{Dattoli}-\cite{Habibullah}. 
In \cite{Johann Cigler},  Cigler introduced 
another family of Hermite polynomials $H_n(x,s)$ generalizing 
the   {\it physicists} and  {\it probabilists} Hermite polynomials as 
\be
\label{samamasaaa}
\sum_{n=0}^\infty\frac{H_n(x,s)}{n!}t^n=e^{xt-s\frac{t^2}{2}}
\ee
with $H_n(x,1)=H_n(x)$ and $H_n(2x,2)={\bf H}_n(x).$ 

In this work, we deal with a construction of  two new families of   $q$ and $(q;p)-$Hermite 
polynomials.

The paper  is organized as follows. In Section \ref{sction1}, we 
give a quick overview on 
the   Hermite polynomials $H_n(x,s)$ introduced in \cite{Johann Cigler}. 
Section  \ref{sction2} is devoted to the construction of
a new family of $q-$Hermite polynomials $H_n(x,s|q)$ 
generalizing the 
discrete $q-$Hermite polynomials. The
inversion formula and relevant properties of these polynomials are computed and discussed.   Their Fourier integral
 transforms are performed in the Section  \ref{sction3}.   
Doubly indexed   Hermite polynomials and   some concluding remarks are introduced in Section  \ref{sction4}.
\section{On the Hermite polynomials $H_n(x,s)$}
\label{sction1}
In \cite{Johann Cigler}, Cigler showed that  the  Hermite polynomials  $H_n(x,s)$ satisfy
\be 
\label{4Hermdiff}
DH_n(x,s)= n\, H_{n-1}(x,s)
\ee
and
the
three term recursion relation 
\be\label{4Hermttr}
H_{n+1}(x,s)=x\, H_n(x,s)-s\,n \,H_{n-1}(x,s),\quad n\geq 1
\ee
with $H_0(x,s):=1.\;D:=d/dx$ is the usual differential operator.  Immediatly, one can see  that
\be
\label{sama:initial}
H_{2n}(0,s)=(-s)^{n}\prod_{k=1}^n(2k-1),\quad H_{2n+1}(0,s)=0.
\ee
The computation of  the first fourth  polynomials gives:
\be 
H_{1}(x,s)= x,\;
H_{2}(x,s)=x^2-s,\;
H_{3}(x,s)=x^3-3\,s\,x,\;
H_{4}(x,s)=x^4-6\,s\,x^2+3\,s^2.
\ee 
More generally, the explicit formula of  $H_{n}(x,s)$ is written  as
\cite{Johann Cigler}
\be 
\label{sasama}
H_{n}(x,s)=n!\sum_{k=0}^{\lfloor\,n/2\,\rfloor}
 \frac{ (-1)^k\,s^k }{ (2k)!! } \frac{ x^{n-2k}  }{(n-2k)!   }
=x^n\,{}_2F_0\Bigg(\begin{array}{c}-\frac{n}{2},\frac{1-n}{2}
\\-\end{array}\Big|\;-\frac{2s}{x^2}\Bigg),
\ee  
where $({}^n_k)=n!/k!(n-k)!$ is a binomial coefficient, 
$n!:=n(n-1)\cdots 2\cdot1,\; (2n)!!:=2n(2n-2)\cdots 2.$\\\\
 The symbol $\lfloor\,x\,\rfloor$ denotes the greatest integer
 in $x$ and ${}_2F_0$ is called the hypergeometric 
series \cite{ASK}.
From 
 (\ref{4Hermdiff}) and  (\ref{4Hermttr}), we have
\be
H_{n }(x,s)= (x-sD)\,H_{n-1 }(x,s),
\ee
where the operator $x-sD $ can  be expressed as
\cite{Johann Cigler}
\be
x-s D = e^{\frac{x^2}{2s}}  (-sD)\,e^{-\frac{x^2}{2s}}.
\ee
The
Rodrigues formula takes the form
\be
\label{sama:abelmarie}
e^{-\frac{x^2}{2s}}\,H_n (x,s)= (-sD)^n\, e^{-\frac{x^2}{2s}}
\ee
while  the  second 
 order
differential equation satisfied by $H_n(x,s)$ is
\be 
\label{second}
\big(sD^2- x D+  n\big)\,H_n(x,s)=0.
\ee
Furthermore, from the  relation (\ref{sasama}) we  derive the result
\be
\label{labelsaj}
H_n(x+sD,s)\cdot(1)=x^n,
\ee
and the  inverse formula for $H_n(x,s)$
\be
\label{labelsaaf}
x^n=n!\,\sum_{k=0}^{\lfloor\,n/2\,\rfloor}
 \frac{s^k   }{  (2k)!! }  \frac{ H_{n-2k}(x,s)   }{(n-2k)! }.
\ee
We then obtain
\be 
\label{props}
\sum_{k,\,n \,(even)}
 \frac{  1    }{(n-k)!\,k! }=\sum_{k,\,n \,(odd)}
 \frac{  1    }{(n-k)!\,k! },\quad 0\leq k\leq n,\quad n\geq 0.
\ee 
From  (\ref{sasama}),  it is also straighforward to note that the polynomials 
$H_n(x,s)$ have an alternative expression given by
\be
\label{mirer}
H_n(x,s)=\exp\left( -s\frac{D^2}{2}\right)\cdot (x^n).
\ee
For any integer $k=0, 1, ..., \lfloor\,n/2\,\rfloor,$ we have the following result
\be 
\label{e7j}
D^{2k}\,H_n(x,s)=\frac{n!}{(n-2k)!} H_{n-2k}(x,s).
\ee
\begin{corollary}
The Hermite polynomials $H_n(x,s)$  obey 
\be
\label{nversion}
\mathcal{T}_n(s, D)\,H_n(x,s)=x^n
\ee
where the polynomial
\be
\label{samapolynom}
\mathcal{T}_n(\alpha,\beta)= \sum_{k=0}^{\lfloor\,n/2\,\rfloor}
 \frac{ 1  }{  (2k)!! } \alpha^k\beta^{2k}.
\ee
\end{corollary}
We are now in a position to formulate and prove the following.
\begin{lemma}
\label{malem}
\be
\mathcal{T}_{2n}(\alpha,\beta)=\frac{(\alpha\beta^2)^n}{(2n)!!}\,{}_2F_0\Bigg(
\begin{array}{c}-n,1\\
-\end{array}\Big|-\frac{2}{\alpha\beta^2}\Bigg)
\ee
and
\be
\mathcal{T}_\infty(\alpha,\beta)=e^{\frac{\alpha\beta^2}{2}}.
\ee
\end{lemma}
{\bf Proof.} From  (\ref{samapolynom}), we have
\bea
\label{calcal}
\mathcal{T}_{2n}(\alpha,\beta)&=&\sum_{k=0}^{n}\frac{1}{(2k)!!}(\alpha\beta^2)^{k}\cr
&=&\frac{(\alpha\beta^2)^n}{(2n)!!}\sum_{k=n}^{\infty}\frac{(2n)!!}{(2k)!!}(\alpha\beta^2)^{k-n}.
\eea
By substituting $m=n-k$ in the latter expression and using various identities, we arrive at
\be 
\label{calcale}
\mathcal{T}_{2n}(\alpha,\beta)
=\frac{(\alpha\beta^2)^n}{(2n)!!}\sum_{m=0}^{\infty}(-n)_m\;\left(\frac{-2}{\alpha\beta^2}\right)^{m},
\ee
where $(a)_j:=a(a+1)\cdots(a+j-1),\; j \geq 1$ and $(a)_0:=1.$
When $n$ goes to $\infty,$ the polynomial (\ref{samapolynom}) takes the form
\be
\mathcal{T}_\infty(\alpha,\beta)= \sum_{k=0}^{\infty}
 \frac{\alpha^k\beta^{2k}  }{  (2k)!! }  = \sum_{k=0}^{\infty}
 \frac{ 1  }{  k!}\left(\frac{\alpha\beta^2}{2}\right)^k
\ee
where $(2k)!!=2^k\, k!$ is used.  $\square$
%

To end this section, let us investigate the Fourier transform of the function $e^{-x^2/2s}H_n(x,s)$. In 
\cite{Johann Cigler}, Cigler has proven that 
\be 
\label{samafouri}
 \frac{1}{\sqrt{2\pi s}}\int_{\mathbb{R}}e^{i x y-\frac{x^2}{2s}} d x=
 e^{-s\frac{y^2}{2}}.
\ee
Hence, 
\be 
\label{samaM}
\frac{1}{\sqrt{2\pi\, s}}\int_{\mathbb{R}}e^{i x y+ i(n- 2k)
 \kappa x- \frac{x^2}{2s}} d x
=e^{-s\frac{y^2}{2}-(n- 2k)s\,y\,\kappa},
\ee 
where   $e^{-2s\kappa^2}=1$.
By  differentiating the relation (\ref{samafouri}) $2n-2k$ times with respect to $y$, one obtains
\be
\frac{1}{\sqrt{2\pi \,s}}\int_{\mathbb{R}}(-1)^{n-k}x^{2n-2k}e^{ix y-\frac{x^2}{2s}} d x=
D^{2n-2k} e^{-s\frac{y^2}{2}}.
\ee
Evaluating the latter expression at $y=0$ and  by making use of  (\ref{sama:abelmarie}), one gets
\be
\label{tate}
\frac{(-1)^{n-k}}{\sqrt{2\pi \,s}}\int_{\mathbb{R}} x^{2n-2k}e^{-\frac{x^2}{2s}} d x=
{D^{2n-2k} e^{-s\frac{y^2}{2}}}_{\big|y=0}={(-s)^{2n-2k}H_{2n-2k}(y,s^{-1})e^{-s\frac{y^2}{2}}}_{\big|y=0}.
\ee
\\
\begin{theorem}
\label{sama:propo1}
The Fourier   transform of the function $e^{-x^2/2 s}H_{n}(x,s )$ is given by 
\be 
\label{samafourier}
  \frac{1 }{\sqrt{2\pi \,s}}\int_{\mathbb{R}}H_{n} (a\, e^{i \kappa x},s)e
^{ix y-\frac{x^2}{2s}}d x=  H_{n} 
 (a\, e^{-s\,\kappa\, y},   s ) e^{-s\frac{y^2}{2}}
\ee
where  $a$ is an arbitrary constant factor. For $y=0,$ we have
\be 
\label{sama:fourier}
  \frac{1 }{\sqrt{2\pi\, s}}\int_{\mathbb{R}}H_{n} (x,s)e
^{-\frac{x^2}{2s}}d x= 0.
\ee
\end{theorem}
{\bf Proof.}
Using (\ref{sasama}) and (\ref{samaM}), we obtain
\bea
 \frac{1}{\sqrt{2\pi\, s}}\int_{\mathbb{R}}H_{n}(a \,e^{i\kappa x},s)
e^{ix y-\frac{x^2}{2\,s}}d x &=&\sum_{k=0}^{\lfloor n/2\rfloor}\frac{(-1)^{k}n!\,s^k\,a^{n-2k}}{(n-2k)!\,(2k)!!} 
\frac{1}{\sqrt{2\pi\, s}}\int_{\mathbb{R}}e^{i x  y+ i(n- 2k)
 \kappa x- \frac{x^2}{2s}} d x \cr
& =& \sum_{k=0}^{\lfloor
n/2\,\rfloor}\frac{(-1)^{k}\,n!\;s^k\,a^{n-2k}}{(n-2k)!\,(2k)!!} 
 e^{-\frac{s}{2}[\kappa (n-2k)+y]^2} \cr   
& =&e^{-s\frac{y^2}{2}} H_{n} (a\, e^{-s\,\kappa\, y}
, s ).  
\eea
 Combining (\ref{sasama}) and (\ref{tate}) for  $n=2n$, we  have
\bea
 \frac{1}{\sqrt{2\pi \,s}}\int_{\mathbb{R}}H_{2n}(x,s)
e^{-\frac{x^2}{2\,s}}d x &=&\sum_{k=0}^n
\frac{(-1)^k\,(2n)!\,s^k}{(2n-2k)!\,(2k)!!}\,
\frac{1}{\sqrt{2\pi\, s}}\int_{\mathbb{R}}x^{2n-2k}\,e^{ixy- \frac{x^2}{2\,s}} d x_{\big|y=0} \cr
&  =&(-1)^n\sum_{k=0}^n
\frac{ (2n)!\,s^k}{(2n-2k)!\,(2k)!!}{ D^{2n-2k} e^{-s\frac{y^2}{2}}}_{\big|y=0} \cr 
&   =&(-1)^{n}\,s^{2n}\,e^{-s\frac{y^2}{2}}\sum_{k=0}^n
\frac{(2n)!\,s^{-k}}{(2n-2k)!\,(2k)!!}{H_{2n-2k}(y,s^{-1})}_{\big|y=0}\cr   
&   =&s^{2n}\,(2n)! \,\sum_{k=0}^n
\frac{(-1)^k}{(2n-2k)!!\,(2k)!!} \cr
&=&0
\eea
where  (\ref{props}) is used. $\square$
\section{New $q-$Hermite polynomials $H_n(x,s|q)$}
\label{sction2}
In this section, we construct 
through the $q-$chain rule 
a new family of $q-$Hermite polynomials  denoted by
$H_n(x,s|q).$
We first introduce some  standard $q-$notations.
For $n\geq 1,\;q\in\mathbb{C}$,
we denote the $q-$deformed number \cite{Ernst} by 
\be
 \{n\}_{q}:=\sum_{k=0}^{n-1}q^k.
\ee
 In the same way, we define the $q-$factorials
\be
 \{n\}_{q}!:=\prod_{k=1}^n \{k\}_{q}, \quad   \{2n\}_{q}!!:=\prod_{k=1}^n \{2k\}_{q},\quad  \{2n-1\}_{q}!!:=\prod_{k=1}^n \{2k-1\}_{q}
\ee
and, by convention, 
\be
  \{0\}_{q}!:=1=: \{0\}_{q}!! \quad \mbox{ and } \quad \{-1\}_{q}!!=1.
\ee 
For any positive number $c,$ the   $q-$Pochhammer symbol $\{c\}_{n,q}$  
is defined as follows:
\be
 \{c\}_{n,q}:=\prod_{k=0}^{n-1} \{c+k\}_{q}
\ee
 while the $q-$binomial coefficients are defined by
\bea
\label{samasama:sa}
{n\atopwithdelims\{\} k}_q:=\frac{ \{n\}_{q}!}{ \{n-k\}_{q}! \{k\}_{q}!}=\frac{(q;q)_n}{(q;q)_{n-k}(q;q)_k}\quad \mbox{ for }\quad
0\leq k\leq n,
\eea
and zero otherwise, where   $(a;q)_n:=\prod_{k=0}^{n-1}(1-a\,q^k),\;(a;q)_0:=1.$ 
\begin{definition}\cite{HahnW, sama}
The Hahn  $q-$addition $\op$ is the function: $\mathbb{C}^3\rightarrow \mathbb{C}^2$
given by:
\be 
 (x,y,q)\mapsto  (x,y)\equiv  x\op y,
\ee 
where 
\bea
\label{addition}
 (x  \op y)^n:&=&(x+y)(x+q \,y)\ldots (x+q^{n-1}\,y)\cr
&=&\sum_{k=0}^n{n
\atopwithdelims\{\} k}_q\, q^{({}^k_2)}\,x^{n-k}\,y^k,\quad n\geq 1,\quad (x\op  y)^0:=1,
\eea
while the $q-$subtraction $\om$ is defined as follows: 
\be 
\label{additionm}
x \om y:=  x\op(-y).
\ee
\end{definition}
Consider a function  $F$ 
\bea
F:D_R  \longrightarrow\mathbb{C},\quad z \longmapsto \sum_{n=0}^\infty c_n\, z^n,
\eea
where $D_R$ is a disc of radius $R.$ We define $F(x\op  y)$ to mean the formal series
\bea
\sum_{n=0}^\infty c_n(x\op  y)^n\equiv\sum_{n=0}^\infty \sum_{k=0}^n\,c_n\,{n
\atopwithdelims\{\} k}_q\, q^{({}^k_2)}\,x^{n-k}\,y^k.
\eea
%
Let $e_q,\;E_q,\; \cos_q$ and $\sin_q$ be the fonctions  defined as follows:
\bea
e_q (x) :&=& \sum_{n=0}^{\infty} \frac{1}{ \{n\}_{q}!}x^n\\
E_q (x) :&=& \sum_{n=0}^{\infty} \frac{q^{n(n-1)/2}}{ \{n\}_{q}!}x^n\cr
\cos_q(x):&=&\frac{e_q (i\,x)+e_q (-i\,x)}{2}= \sum_{n=0}^{\infty} \frac{(-1)^n}{ \{2n\}_{q}!}x^{2n}\\
\sin_q(x):&=&\frac{e_q (i\,x)-e_q (-i\,x)}{2\,i}= \sum_{n=0}^{\infty} \frac{(-1)^n}{ \{2n+1\}_{q}!}x^{2n+1}.
\eea
We immediately obtain the following rules for the product of two exponential functions
\be
e_q (x) E_q (y) = e_q ( x \op y ).
\ee
The new family of $q-$Hermite polynomials $H_n(x,s|q) $ can be determined 
by the  generating function
\be
\label{functiongeneratrice}
e_q\big(  tx \ominus_{q,q^2} st^2/ \{2\}_{q} \big) =
 e_q(  tx)E_{q^2} (- st^2/ \{2\}_{q} ) := \sum_{n=0}^{\infty} \frac{H_n(x,s|q) }{ \{n\}_{q}!}t^n,\quad |t| <1,
\ee
where \cite{sama}
\be
(a\ominus _{q, q^{2}}b )^n := \sum_{k=0}^{n}   \frac{ \{n\}_{q}!}{ \{n-k\}_{q}! \, \{k\}_{q^{2} }!}
(-1)^k q^{k(k-1)} a^{n-k} \,b^{k},\quad (a\ominus_{q, q^{2}} b )^0  := 1.
\ee
Performing the $q-$derivative $D_x^q$ of   both sides 
of  (\ref{functiongeneratrice}) with respect to $x$, one obtains
\bea
\label{lowering}
D_x^q\, H_n(x,s|q)= \{n\}_{q}\, H_{n-1}(x,s|q),
\eea
where
\be
\label{lowm}
D_x^q\, f(x)=\frac{f(x)-f(qx)}{(1-q)x}
\ee
 satisfying
\bea
D_x^q (a\, x \op  b)^n= \{n\}_{q}\,(a\, x \op  b)^{n-1}.
\eea

Recall \cite{Ismails}  that the Al-Salam-Chihara polynomials $P_n(x;a,b,c)$ satisfy the following  recursion relation:
\be
\label{sama:al}
P_{n+1}(x;a,b,c)=(x-a\,q^n)\,P_n(x;a,b,c)-(c+b\,q^{n-1})\, \{n\}_{q}\,P_{n-1}(x;a,b,c)
\ee
with  $P_{-1}(x;a,b,c)=0$  and $P_0(x;a,b,c)=1$.\\
Performing the $q-$derivative 
of both sides of
  (\ref{functiongeneratrice}) with respect to $t$, we have
\bea
\label{raising}
H_{n+1}(x,s|q)=  x\, H_{n}(x,s|q)-s\,  \{n\}_{q} \,q^{n-1}\,H_{n-1}(x,s|q),\quad n\geq 1
\eea
with $H_0(x,s|q):=1.$ \\
By setting
$a=0=c$ and $b=s$ in (\ref{sama:al}), one obtains the   recursion relation (\ref{raising}).
 From the latter equation,
one can see  that
\be
\label{sama:qinitial}
H_{2n}(0,s|q)=(-s)^n\,q^{n(n-1)}\, \{2n-1\}_{q}!!,\quad  H_{2n+1}(0,s|q)=0.
\ee
The first fourth new   polynomials are given by
\bea
H_1(x,s|q) &=&  x,\\
H_2(x,s|q)& =&  x^2 -s,\\
H_3(x,s|q)&=&  x^3 - \{3\}_{q}s x,\\
H_4(x,s|q)&=&  x^4 -(1+q^2) \{3\}_{q}sx^2+ q^2\,\{3\}_{q}s^2.
\eea 
More generally, we have the following.
\begin{theorem}
 The explicit formula for the new   Hermite polynomials $H_n(x,s|q)$ is given by
\bea
\label{sama:qhermite}
H_n(x,s|q) &=& \sum_{k=0}^{\lfloor\,n/2\,\rfloor}
 \frac{ (-1)^k q^{k(k-1)}  \{n\}_{q}!  }{ \{n-2k\}_{q}! \, \{2k\}_{q}!! }s^k x^{n-2k}\\
&=&x^n\,{}_2\phi_0\Bigg(\begin{array}{c}q^{-n},q^{1-n}\\
-\end{array}\Big|\;q^2;\;\frac{sq^{2n-1}}{(1-q)x^2}\Bigg),
\eea
where ${}_2\phi_0$ is the $q-$hypergeometric series \cite{ASK}.
\end{theorem}
{\bf Proof.} Expanding the generation function 
given in (\ref{functiongeneratrice}) in Maclaurin series, we have
\bea
\label{sammas}
 e_q( t\,x)E_{q^2} (- s\,t^2 /\{2\}_{q}) &=& \sum_{k=0}^{\infty} \frac{( x\,t)^k }{\{k\}_{q}!} 
\sum_{m=0}^{\infty} \frac{(-1)^m q^{m(m-1)} }{\{m\}_{q^2} ! } \left(\frac{s\,t^2}{\{2\}_{q}}\right)^m \cr
 &=& \sum_{k=0}^{\infty} \sum_{m=0}^{\infty}  \frac{(-1)^m q^{m(m-1)} 
 x^k }{\{k\}_{q}! \,\{m\}_{q^2} ! } \left(\frac{s}{\{2\}_{q}}\right)^m t^{k+ 2m}.
\eea
By substituting 
\be
k +2m =n ~\Rightarrow~ m \le\lfloor\,n/2\,\rfloor,
\ee
and 
\be
\{2\}_{q}\,\{m\}_{q^2}=\{2m\}_{q}
\ee
in (\ref{sammas}),   we have
\be
 e_q( t\,x)E_{q^2} (- s\,t^2/\{2\}_{q} ) = \sum_{n=0}^{\infty} \lb
 \sum_{m=0}^{\lfloor\,n/2\,\rfloor}   \frac{(-1)^m q^{m(m-1)}
 s^m\,x^{n-2m} }{\{n-2m\}_{q}!\, \{2m\}_{q} !! }\rb  t^{n},
\ee
which achieves the proof. $\square$

In the limit case when $x\to \{2\}_{q}\,x,\;s\to (1-q)\,\{2\}_{q},$ the
 polynomials $H_n(x,s|q)$ are reduced to $H_n^q(x)$ investigated by Chung et {\it al} \cite{sama}. When $s\to 1-q$, they are 
reduced to   the discrete $q-$Hermite I polynomials \cite{ASK}.\\
The  relation (\ref{raising}) allows us to write
\be
H_{n}(x,s|q)=(x-sq^N\circ D_x^q)\,H_{n-1}(x,s|q),
\ee
where
 the operator $N$ acts on the polynomials $H_n(x,s|q)$ as follows:
\be 
NH_n(x,s|q):= n\,H_n(x,s|q),\quad q^N\circ D_x^q=D_x^q\circ q^{N-1}.
\ee
It is straightforward to show that the polynomials (\ref{sama:qhermite}) satisfy the 
following $q-$difference equation
\be
\label{sama:hounk}
\big(s\, (D_x^q)^2 -x\,q^{2-n}\,D_x^q+q^{2-n}\,\{n\}_{q}\big)\,H_n(x,s|q)=0.
\ee
In the limit case when $q$ goes to $1$, the $q-$difference 
equation (\ref{sama:hounk})   reduces to the well-known differential equation (\ref{second}). 
For $n$ even or odd, the  polynomials $H_{n}(x,s|q)$ obey   the following  generating functions
\be
\sum_{n=0}^\infty\frac{H_{2n}(x,s|q)}{\{2n\}_{q}!}(-t)^n=\cos_q(x\sqrt{t})\,E_{q^2}(s\,t/\{2\}_{q}),\; |t| <1
\ee
or
\be
\sum_{n=0}^\infty\frac{H_{2n+1}(x,s|q)}{\{2n+1\}_{q}!}(-t)^n=
\frac{1}{\sqrt{t}}\sin_q(x\sqrt{t})\,E_{q^2}(s\,t/\{2\}_{q}), \; |t| <1,
\ee
respectively.
\begin{theorem}
\label{theo:poly}
The polynomials $H_n(x,s|q)$ can be expressed as
\be
H_n(x,s|q) =\prod_{k=1}^{n} \big( x -s\,q^{n-1-k}\, D_x^q \big)\cdot(1).
\ee
We also have
\bea
\label{sama:reso}
H_n \left(x+sq^N\circ D_x^q,s|q\right)\cdot(1)=x^n.
\eea
\end{theorem}
{\bf Proof.}  
Since (\ref{lowering}) and (\ref{raising}) are satisfied,  we have
\bea
H_{n}(x,s|q)&=&  x\ H_{n-1}(x,s|q)-s\, q^{n-2}\, \{n-1\}_{q}\, H_{n-2}(x,s|q)\cr
&=& x \,H_{n-1}(x,s|q)-s  q^{n-2}\,D_x^q\,H_{n-1}(x,s|q).
\eea
The rest holds by induction on $n$.\\
To prove the relation   (\ref{sama:reso}) we  replace $x^{n-2k}$ in (\ref{sama:qhermite}) by 
$(x+sq^N\circ D_x^q)^{n-2k}$ and apply the   corresponding linear operator
to $1$. The relation (\ref{sama:reso}) is true for $n=0$ and $n=1.$
For $n=2$, we have
\bea
H_2\left(x+sq^N\circ D_x^q,s|q\right)\cdot(1)&=&\left(x+sq^N\circ D_x^q\right)^2\cdot(1)-s\cr
&=&\left(x+sq^N\circ D_x^q\right)\cdot (x)-s\cr
&=&x^2.
\eea
  Assume
that  (\ref{sama:reso}) is true for $n-1,\;n\geq 3.$
 Then we must prove that 
\be
H_{n} \left(x+sq^N\circ D_x^q,s|q\right) \cdot(1)=x^{n}.
\ee
From 
 (\ref{raising}),  we have
\bea
H_{n} \left(x+sq^N\circ D_x^q,s|q\right) \cdot1&=& \left(x+sq^N\circ D_x^q\right)H_{n-1}
 \left(x+sq^N\circ D_x^q,s|q\right) \cdot(1)\cr
&-&s\{n-1\}_{q}\,q^{n-2}H_{n-2}  \left(x+sq^N\circ D_x^q,s|q\right)\cdot(1)\cr
&=&\left(x+sq^N\circ D_x^q\right)\cdot x^{n-1}
-s\{n-1\}_{q}\,q^{n-2} x^{n-2}\cr
&=&x^{n}
\eea
which achieves the proof. $\square$\\
From the   {\bf Theorem} \ref{theo:poly}, we obtain the following.
\begin{corollary}
\label{sama:theorem}
The   polynomials (\ref{sama:qhermite}) have the following inversion formula
\be
\label{inversion}
x^n = \{n\}_{q}!\sum_{k=0}^{\lfloor\,n/2\,\rfloor}
 \frac{ q^{k(k-1)} \;s^k  }{  \{2k\}_{q}!! }\frac{ H_{n-2k}(x,s|q) }{\{n-2k\}_{q}!  }.
\ee
\end{corollary}
{\bf Proof.} Let  $h_n^q(x,s)$ be the polynomial defined by
\be
h_n^q(x,s)=\left(x+sq^N\circ D_x^q\right)^n\cdot(1).
\ee
Note that
$
h_n^q(x,-s)=H_n(x,s|q).
$
From (\ref{sama:reso}), we have
\bea
x^n&=&\sum_{k=0}^{\lfloor\,n/2\,\rfloor}
 \frac{(-1)^k q^{k(k-1)}\{n\}_{q}!  }{\{n-2k\}_{q}! \{2k\}_{q}!! }s^k 
\left(x+sq^N\circ D_x^q\right)^{n-2k}\cdot(1)\cr
&=&\sum_{k=0}^{\lfloor\,n/2\,\rfloor}
 \frac{ q^{k(k-1)}\{n\}_{q}! \,s^k }{\{n-2k\}_{q}!\{2k\}_{q}!! } h_{n-2k}^q(x,-s)
\eea
which achives the proof.
$\square$

From (\ref{lowering}), one readily  deduces  that, for integer 
powers $k=0, 1, ..., \lfloor\,n/2\,\rfloor$ of the operator 
$D_x^q,$
\be 
\label{e7}
(D_x^q)^{2k}H_n(x,s|q)=\gamma_{n,k}(q)H_{n-2k}(x,s|q),\quad 
\gamma_{n,k}(q)=\frac{\{n\}_{q}!}{\{n-2k\}_{q}!}.
\ee
Therefore, we have the following decomposition of  unity 
\be
\sum_{k=0}^{\lfloor\,n/2\,\rfloor}\frac{  (-1)^kq^{k(k-1)} s^k }{  \{2k\}_{q}!! 
}(D_x^q)^{2k}\sum_{m=0}^{\lfloor\,n/2\rfloor}
   \frac{ q^{m(m-1)}   s^m }{   \{2m\}_{q}!! }(D_x^q)^{2m}={\bf 1}
\ee
and the new 
  $q-$Hermite polynomials $H_n(x,s|q)$ 
obey 
\bea
\label{rinversion}
\mathcal{L}_n(s, D_x^q|q)H_n(x,s|q)=x^n 
\eea
where the polynomial $\mathcal{L}_n(\alpha,\beta|q)$ is defined as follows:
\be
\label{sama:polynom}
\mathcal{L}_n(\alpha,\beta|q)= \sum_{k=0}^{\lfloor\,n/2\,\rfloor}
 \frac{ q^{k(k-1)}   }{  \{2k\}_{q}!! } \alpha^k\beta^{2k}.
\ee
This  polynomial is   essentially the $(q,n)-$exponential 
function  $E_{q,n}(x)$ investigated by Ernst \cite{Ernst}, i.e.,
$\mathcal{L}_{n-1}(\alpha,\beta|q)=E_{q^{-2},\lfloor\,n/2\,\rfloor}(\alpha\beta^2/\{2\}_{q}).$ 
We are now in a position to formulate and prove the following.
\begin{lemma}
From the polynomial  (\ref{sama:polynom})  we have
\be
\mathcal{L}_{2n}(\alpha,\beta|q)=\frac{(\alpha\beta^2)^nq^{n(n-1)}}{\{2n\}_{q}!!}{}_3\phi_2\Bigg(
\begin{array}{c}q^{-n},-q^{-n},q\\
0,0\end{array}\Big|q;-\frac{q^2}{(1-q)\alpha\beta^2}\Bigg)
\ee
and
\be
\mathcal{L}_\infty(\alpha,\beta|q)=E_{q^2}(\alpha\beta^2/\{2\}_{q}).
\ee
\end{lemma}
{\bf Proof.} As it is  defined in (\ref{sama:polynom}),    we have
\begin{eqnarray}
\mathcal{L}_{2n}(\alpha,\beta|q)&=&\sum_{k=0}^{n}\frac{q^{ k(k-1)}}{\{2k\}_{q}!!}(\alpha\beta^2)^{k}\cr
&=&\frac{(\alpha\beta^2)^n}{\{2n\}_{q}!!}\sum_{k=n}^{\infty}\frac{q^{ k(k-1)}\{2n\}_{q}!!}{\{2k\}_{q}!!}(\alpha\beta^2)^{k-n}.
\end{eqnarray}
By substituting
$m=n-k$ in the latter expression, we arrive at
\begin{eqnarray}
\mathcal{L}_{2n}(\alpha,\beta|q)
&=&\frac{(\alpha\beta^2)^n}{\{2n\}_{q}!!}\sum_{m=0}^{\infty}
\frac{q^{ (n -m)(n -m-1)}\{2n\}_{q}!!}{\{2n-2m\}_{q}!!}(\alpha\beta^2)^{-m}\cr
&=&\frac{(\alpha\beta^2)^nq^{n(n-1)}}{\{2n\}_{q}!!}\sum_{m=0}^{\infty}
 (q^{-2n};q^2)_m\left(-\frac{q^2}{(1-q)\alpha\beta^2}\right)^m.
\end{eqnarray}
When $n\to\infty$,   (\ref{sama:polynom}) takes the form
\be
\mathcal{L}_\infty(\alpha,\beta|q)= \sum_{k=0}^{\infty}
 \frac{ q^{k(k-1)}   }{ \{2n\}_{q}!! } (\alpha\beta^2)^k = \sum_{k=0}^{\infty}
 \frac{ q^{k(k-1)}   }{  \{k\}_{q^2 }!}\left(\frac{\alpha\beta^2}{\{2\}_{q}}\right)^k
\ee
which achieves the proof. $\square$

In the limit, when $q\to1$, the polynomial  $\mathcal{L}_n(\alpha,\beta|q)$ is reduced 
to the classical one's  $\mathcal{T}_n(\alpha,\beta)$, i.e., $\lim_{q\to 1}
\mathcal{L}_n(\alpha,\beta|q)=\mathcal{T}_n(\alpha,\beta),\;\forall\,n.$

\section{Fourier transforms of the new $q-$Hermite 
polynomials $H_n(x,s|q)$ }
\label{sction3}
In this section, 
we compute the Fourier integral
transforms associated to the new $q-$Hermite 
polynomials $H_n(x,s|q)$. 
\subsection{ $q^{-1}-$Hermite polynomials $H_n(x,s|q^{-1})$ }
Let us rewrite the  new $q-$Hermite polynomials  
(\ref{sama:qhermite}) in the following form
\be 
 \label{sama:pipp}
H_{n}(x,s|q)= \sum_{k=0}^{\lfloor\,n/2\,\rfloor}
c_{n,k}(q)\,s^{k}x^{n-2k},               
 \ee
where  the associated  coefficients 
$c_{n,\,k}(q)$ are given by
\be 
 \label{sama:coef}
c_{n,\,k}(q):= \frac{ (-1)^k q^{k(k-1)}\,\{n\}_{q}! 
 }{\{n-2k\}_{q}! \,\{2k\}_{q}!! }.
\ee
By a direct computation, one can easily check that these coefficients 
satisfy the following recursion relation
\be 
\label{rett}
c_{n+1,\,k}(q)=c_{n,\,k}(q)-q^{n-1}\{n\}_{q}\,c_{n-1,\,k-1}(q),
\ee
with $
c_{0,\,k}(q)=\delta_{0,k}, \;\; c_{n,\,0}(q)=1.$\\
From the definition of the $q-$binomial coefficients
 in (\ref{samasama:sa}),  it is not hard to derive an
inversion formula
\bea
\label{samasamasa}
{n\atopwithdelims\{\} 2k}_{q^{-1}}=
q^{2k(2k-n)}\,{n\atopwithdelims\{\} 2k}_q,\qquad
0\leq k\leq \lfloor n/2\rfloor.
\eea
Then, one readily deduces that
\be 
\label{sama:labas}
c_{n,\,k}(q^{-1})= q^{k(k+3-2n)}c_{n,\,k}(q),
\ee
allowing to define the $q^{-1}-$Hermite 
polynomials $H_n(x,s|q^{-1})$ in the following form
\be 
 \label{sama:inm}
H_n(x,s|q^{-1}):= \sum_{k=0}^{\lfloor\,n/2\,\rfloor}\,c_{n,\,k}(q^{-1})\,s^{k}x^{n-2k}.                      \ee
The  recursion relation 
\be 
\label{sama:lab}
c_{n+1,\,k}(q^{-1})= q^{-2k}\,c_{n,\,k}(q^{-1})- q^{3-n-2k}\,\{n\}_{q}\,c_{n-1,\,k-1}(q^{-1}),\quad n\geq 1
\ee
is valid for the coefficients (\ref{sama:labas}) 
with $c_{0,\,k}(q^{-1})=q^{k(k+3)}\,\delta_{0,k}, \;\; c_{n,\,0}(q^{-1})=1.$ \\
Since  (\ref{sama:lab}) is satisfied, the $q^{-1}-$Hermite polynomials $H_n(x,s|q^{-1} )$ obey 
the relation
\be 
 \label{samanm}
H_{n+1}(x,s|q^{-1})=xH_n(x,sq^{-2}|q^{-1})-sq^{1-n}\{n\}_{q}\,H_{n-1}(x,sq^{-2}|q^{-1}),\quad n\geq 1,
\ee
with $H_0(x,sq^{-2}|q^{-1}):=1.$\\
The action of the operator $D_x^q$ on 
the   polynomials (\ref{sama:inm}) is given by
\bea
 \label{samanma}
D_x^q\,H_n(x,s|q^{-1})=\{n\}_{q}\,H_{n-1}(x,sq^{-2}|q^{-1}).
\eea
Let $\epsilon$ 
denote the operator which maps $f(s)$ to $f(qs)$. Then, from (\ref{samanm}) and (\ref{samanma}) one 
can establish that
\be
H_{n}(x,s|q^{-1}) =\prod_{k=1}^{n} \big(x\,\epsilon^{-2}-sq^{k+1-n}D_x^q\big)\cdot(1).
\ee
\subsection{Fourier transforms of the new $q-$Hermite polynomials $H_n(x,s|q)$ }
Considering
the well-known Fourier transforms (\ref{samafouri})
for the Gauss exponential function $e^{-x^2/2s},$  the 
Fourier integral transforms
for the exponential function $\exp(i(n-2k)\kappa x-x^2/2s)$  is computed as follows:
\be 
\frac{1}{\sqrt{2\pi s}}\int_{\mathbb{R}}e^{i x  y+ i(n- 2k)
 \kappa x- \frac{x^2}{2s}} d x
=q^{\frac{n^2}{4}+k(k-n)}e^{-s\frac{y^2}{2}-(n- 2k)sy\kappa},
\ee 
where   $q=e^{-2s\kappa^2} \leq1$ and $0\leq \kappa < \infty$.
\begin{theorem}
\label{sama:propo1}
The new $q-$Hermite polynomials $H_n(x,s|q)$ and 
$H_n(x,s|q^{-1})$ defined in (\ref{sama:pipp}) and (\ref{sama:inm}),
respectively, are connected  by the integral Fourier transform  of the
following form
\be 
\label{sama:fourier}
  \frac{1 }{\sqrt{2\pi s}}\int_{\mathbb{R}}H_{n} (b e^{i \kappa x},s|q)e
^{ix y-\frac{x^2}{2s}}d x= q^{\frac{n^2}{4}}H_{n}
 (b e^{-s\kappa y}, q^{n-3} s| q^{-1})\, e^{-s\frac{y^2}{2}}
\ee
where $b$ is an arbitrary constant factor.
\end{theorem}
 {\bf Proof.} To prove this theorem, let us make use  of (\ref{sama:pipp}) and 
evaluate the left hand side of (\ref{sama:fourier}). Then,
\bea
 \frac{1}{\sqrt{2\pi s}}\int_{\mathbb{R}}H_{n}(b e^{i\kappa x},s|q)
e^{ix y-\frac{x^2}{2s}}d x &=&\sum_{k=0}^{\lfloor n/2\rfloor}c_{n,k} (q)s^{k}
b^{n-2k}\,\frac{1}{\sqrt{2\pi s}}\int_{\mathbb{R}}e^{i x  y+ i(n- 2k)
 \kappa x- \frac{x^2}{2s}} d x \cr
&  =&\sum_{k=0}^{\lfloor
n/2\,\rfloor}c_{n,k}(q)s^k b^{n-2k}e^{-\frac{s}{2}[\kappa (n-2k)+y]^2} \cr
&  = &q^{\frac{n^2}{4}}\sum_{k=0}^{\lfloor
n/2\,\rfloor}c_{n,k}(q)q^{-k(n-k)}s^k b^{n-2k}e^{-s\frac{y^2}{2}-(n- 2k)sy\kappa} \cr  
&  = &q^{\frac{n^2}{4}}\sum_{k=0}^{\lfloor
n/2\,\rfloor}c_{n,k}(q^{-1})(q^{n-3}s)^k (be^{-sy\kappa})^{n-2k} e^{-s\frac{y^2}{2}}\cr             
& =& q^{\frac{n^2}{4}}H_{n}
 (b e^{-s\kappa y}, q^{n-3} s |q^{-1})\, e^{-s\frac{y^2}{2}}.  
\eea
$\square$

\section{Doubly indexed  Hermite polynomials $\mathcal{H}_{n,p}(x,s|q)$}
\label{sction4}
In this section, we construct a novel family  of   Hermite polynomials  
  called  {\it doubly indexed  Hermite polynomials}, $\mathcal{H}_{n,p}(x,s|q).$  
First, let us 
defined the 
$(q;p)-$shifted factorials $ (a;q)_{pk}$
and  the $(q;p)-$number as follows:
\be
\label{sama:newdef}
(a;q)_0:=1,\quad (a;q)_{pk}:=(a,aq,\cdots,aq^{p-1};q^p)_k, \quad p\geq 1,\;k=1,2,3,\cdots
\ee
and 
\be
 \{pk\}_{q}:=\frac{1-q^{pk}}{1-q},\quad    \{pk\}_{q}!!:=\prod_{l=1}^k  \{pl\}_{q},\quad  \{0\}_{q}!!:=1,
\ee
respectively.
\begin{definition}
For a positive integer $p$,   a class of doubly indexed  
Hermite polynomials $\big\{\mathcal{H}_{n,p}\big\}_{n,p}$  is defined 
such that
\be
\label{sama:doublyH}
\mathcal{H}_{n,p}(x,s|q) := E_{q^p} \lb - s\frac{(D_x^q )^p}{ \{p\}_{q} }  \rb\cdot (x ^n).
\ee
\end{definition}
If  $p=2,$ a subclass of the polynomials 
(\ref{sama:doublyH}) is reduced to the class of polynomials 
  (\ref{sama:qhermite}). 
More generally, 
 their explicit formula  is given by
\bea
\label{sama:qdoublyH}
\mathcal{H}_{n,p}(x,s|q)&=&   \{n\}_{q}!\sum_{k=0}^{\lfloor\,n/p\,\rfloor}
 \frac{ (-1)^k q^{p({}^k_2)}  \,s^k}{  \{pk\}_{q}!!}\frac{  x^{n-pk}  }{  \{n-pk\}_{q}! }\\
&=&x^n{}_p\phi_0\left(\begin{array}{c}q^{-n},q^{-n+1},\cdots,q^{-n+p-1}\\
-\end{array}\Big|\;q^p;\;\frac{sq^{p(n+(1-p)/2)}}{(1-q)^{p-1}x^p}\right),
\eea
where ${}_p\phi_0$ is the $q-$hypergeometric series \cite{ASK}.

Since $D_x^qe_q(\omega\,x)=\omega\,e_q(\omega\,x),$  we  derive   
the generating function of the 
polynomials (\ref{sama:doublyH})
as 
 \be
\label{sama:functionqdoublygeneratrice}
 f_q(x,s;p):=e_q(  tx)E_{q^p} (- st^p/ \{p\}_{q} ) =
 \sum_{n=0}^{\infty} \frac{\mathcal{H}_{n,p}(x,s|q) }{ \{n\}_{q}!}t^n, \quad |t| < 1.
\ee
These  polynomials   are the solutions of the $q-$analogue of the generalized  heat equation \cite{Dattoli}
\be
(D_x^q)^p f_q(x,s;p)= -\{p\}_qD_s^qf_q(x,s;p),\quad  f_q(x,0;p)=x^n.
\ee
For any real number $c$ and a positive integer $p$, $|q| <1,$ we have
\be
  \sum_{n=0}^{\infty} \frac{\{c\}_{n,q}\mathcal{H}_{n,p}(x,s|q) }{\{n\}_{q}!}t^n=
\frac{1}{(xt;q)_c}{}_p\phi_p\Bigg(\begin{array}{c}q^{c},q^{c+1},\cdots,q^{c+p-1}\\
xtq^{c},xtq^{c+1},\cdots,xtq^{c+p-1}\end{array}\Big|\;q^p;\; \frac{s\,t^p}{(1-q)^{p-1}}\Bigg), \;|xt| <1.
\ee
Performing the $q-$derivative   of   both sides 
of  (\ref{sama:functionqdoublygeneratrice}) with respect to $x$ and $t$, one obtains
\be 
\label{sa:lowering}
D_x^q \,\mathcal{H}_{n,p}(x,s|q)= \{n\}_{q} \,\mathcal{H}_{n-1,p}(x,s|q)
\ee 
and
\be 
\label{sa:raising}
\mathcal{H}_{n+1,p}(x,s|q)=  x \mathcal{H}_{n,p}(x,s|q)-s q^{n-p+1} 
\{n\}_{q} \{n-1\}_{q} \cdots \{n-p+2\}_{q}  \mathcal{H}_{n-p+1}(x,s|q),\; n\geq 1,
\ee 
  with $\mathcal{H}_{0,p}(x,s|q):=1.$ 
 The  polynomials  (\ref{sama:doublyH}) obey   the following $p-$th order difference equation
\be
\Big( s\,(D_x^q)^p-q^{p-n}\,x\,D_x^q+q^{p-n}\, \{n\}_{q}\Big)\,\mathcal{H}_{n,p}(x,s|q)=0.
\ee

\section{Concluding remarks}
 In this paper, we   have constructed a  family of new $q-$Hermite polynomials 
$H_n(x,s|q)$.
Several  properties related to these 
polynomials have been computed and discussed.
Finally, we have constructed
a novel family of Hermite polynomials $\mathcal{H}_{n,p}(x,s|q)$ 
called doubly indexed Hermite polynomials. 

In the limit cases, when $q$  goes to $1$ and $s$ goes to $-py,$ the 
polynomials $\mathcal{H}_{n,p}(x,s|q)$ are reduced to
the higher-order Hermite polynomials, sometimes
called the Kamp\'e de F\'eriet or the Gould Hopper polynomials \cite{Dattoli}-\cite{GOULDHOPPER}, i.e.,
\be
\mathcal{H}_{n,p}(x,-py|1)\equiv g_{n}^p(x,y):=   n!\sum_{k=0}^{\lfloor\,n/p\,\rfloor}
 \frac{  y^kx^{n-pk}}{  k! \, (n-pk)! }.
\ee
 When  $q$  goes to $1,\;x\to px$ and $s\to p,$ the polynomials $\mathcal{H}_{n,p}(x,s|q)$
become   the Hermite 
polynomials  investigated by Habibullah 
and  Shakoor \cite{Habibullah}, i.e.,
\be
\mathcal{H}_{n,p}(px,p|1)\equiv S_{p,n}(x):=   n!\sum_{k=0}^{\lfloor\,n/p\,\rfloor}
 \frac{  (-1)^k(px)^{n-pk}}{  k! \, (n-pk)! }.
\ee 
For $p=2$, the doubly indexed polynomials   $\mathcal{H}_{n,p}(x,s|q)$ are reduced to the new $q-$Hermite polynomials $H_{n}(x,s|q),$ i.e.,
$\mathcal{H}_{n,2}(x,s|q)\equiv H_{n}(x,s|q)$. 

\section*{Acknowledgements}
This work is partially supported by the Abdus Salam International
Centre for Theoretical Physics (ICTP, Trieste, Italy) through the
Office of External Activities (OEA)-\mbox{Prj-15}. The ICMPA
is in partnership with
the Daniel Iagolnitzer Foundation (DIF), France.

\end{document}